\documentclass{amsart}
\usepackage{epsfig, times}
\usepackage{mathrsfs}

\input epsf
\def\relabelbox{%
  \hbox\bgroup%
}%
\def\endrelabelbox{%
\egroup%
}%
\def\relabel #1#2 {%
  \special{ps:/a {} def}%
  \smash{\rlap{#2}}%
}%
\def\adjustrelabel <#1,#2> #3#4 {%
  \special{ps:/a {} def}%
  \smash{\rlap{\kern #1 \raise #2\hbox{#4}}}%
}%
\def\extralabel <#1,#2> #3 {\smash{\rlap{\kern #1 \raise #2\hbox{#3}}}}%

\begin{document}

\newtheorem{thm}{Theorem}[section]
\newtheorem{lem}[thm]{Lemma}
\newtheorem{slem}[thm]{Sublemma}
\newtheorem*{main_thm}{Theorem {\ref{coarse_package}}}
\newtheorem{prop}[thm]{Proposition}

\newtheorem{cor}[thm]{Corollary}
\newtheorem{conj}[thm]{Conjecture}
\newtheorem{qn}[thm]{Question}
\newtheorem{claim}[thm]{Claim}
\newtheorem{prob}[thm]{Problem}
\newtheorem{iprob}{Problem}
\theoremstyle{definition}
\newtheorem{defn}[thm]{Definition}

\theoremstyle{remark}
\newtheorem{rmk}[thm]{Remark}
\newtheorem{exa}[thm]{Example}

\def\square{\hfill${\vcenter{\vbox{\hrule height.4pt \hbox{\vrule width.4pt
height7pt \kern7pt \vrule width.4pt} \hrule height.4pt}}}$}

\def\R{\mathbb R}
\def\Z{\mathbb Z}
\def\O{\mathscr O}
\def\P{\mathscr P}
\def\F{\mathscr F}
\def\G{\mathscr G}
\def\C{\mathscr C}
\def\M{\mathscr M}
\def\L{\mathscr L}
\def\E{\mathbb E}
\def\e{\mathscr E}
\def\QP{\mathbb{QP}}
\def\RP{\mathbb{RP}}
\def\H{\mathbb{H}}
\def\til{\widetilde}
\def\im{\text{im}}
\def\homeo{\text{Homeo}}
\def\Hom{\text{Hom}}
\def\u{{\text{univ}}}

\newenvironment{pf}{{\it Proof:}\quad}{\square \vskip 12pt}
\newenvironment{spf}{{\it Sketch of proof:}\quad}{\square \vskip 12pt}
\newenvironment{pmthm}{{\it Proof of theorem~\ref{coarse_package}:}\quad}{\square \vskip 12pt}

\title{Irrational Stable Commutator Length in Finitely Presented Groups}

\author{Dongping Zhuang}
\address{Department of Mathematics \\ Caltech \\ Pasadena CA 91125}
\email{dongping@caltech.edu}
\maketitle

\begin{abstract}
We give examples of finitely presented groups containing elements
with irrational (in fact, transcendental) stable commutator length,
thus answering in the negative a question of M. Gromov. Our examples
come from 1-dimensional dynamics, and are related to the generalized
Thompson groups studied by M. Stein, I. Liousse and others.
\end{abstract}

\section{Introduction}
Let $G$ be a group, and let $a$ be an element of the commutator
subgroup, which we always denote by $[G,G]$. The {\em commutator
length} of $a$, which we denote $cl(a)$, is defined to be the
minimum number of commutators whose product is equal to $a$. That
is,
$$cl(a) = \min \lbrace n | a = [b_1,c_1] \cdots[b_n,c_n] \rbrace$$
$cl(a)$ is a subadditive function, so the limit of $\frac{cl(a^n)}
n$, as $n \to \infty$, exists.
\begin{defn}
$G, a \in [G, G]$, the {\em stable commutator length} of $a$,
denoted by $scl(a)$, is defined to be
$$scl(a) = \lim_{n \to \infty} \frac {cl(a^n)} n$$
Set $scl(a)=\infty$, if no power of $a$ is in $[G,G]$.
\end{defn}\
Commutator length and stable commutator length in groups have long
been studied, often under the name of the {\em genus} problem. If
$G=\pi_1(V)$ for some aspherical space $V$, and $\gamma$ is a loop
in $V$ representing the conjugacy class of an element $a$, then
the commutator length of $a$ is the minimal genus of a surface $S$
for which there is a map $f:S \to V$ taking $\partial S$ to
$\gamma$. More generally, given a class $H \in H_2(V,\gamma)$, one
can ask for the least genus immersed surface in the relative class
$H$. Stabilizing, one obtains a norm on $H_2(V,\gamma)$. If $G$ is
finitely presented, $H_2(V,\gamma)$ is finitely dimensional, and
one can try to minimize this stable norm on the subspace which is
the preimage of the class of $[\gamma]$ under the boundary map
$H_2(V,\gamma) \to H_1(\gamma)$. This infimum is the stable
commutator length of $a$. \

M. Gromov (in \cite{gr2} $6.C_2$) asked the question of whether such
a stable norm on $H_2(V,\gamma)$, or in our content, the stable
commutator length in a finitely presented group, is always rational,
or more generally, algebraic. The purpose of this note is to give
simple (and even natural) examples which show that the stable
commutator length in finitely presented groups can be
transcendental.\

We now state the contents of this note. In $\S 2$ we state the
fundamental {\em Duality Theorem} of C. Bavard \cite{ba}, which
gives a precise relationship between stable commutator length and
homogeneous quasimorphisms on groups. In $\S 3$ we give our examples
and demonstrate that they have the desired properties, which are
based on the work of M. Stein \cite{st}, D. Calegari \cite{ca3} and
I. Liousse \cite{li}. The examples are central extensions of certain
(finitely generated) groups of piecewise linear homeomorphisms of
the circle. C. Bavard's duality theorem connects dynamics (rotation
numbers, as studied by I. Liousse \cite{li}) with stable commutator
length. We include an appendix presenting basic properties of
rotation numbers. They are used in $\S 2$ and $\S 3$.\

\subsection{Acknowledgements}
The fact that these examples have irrational stable commutator
length was first conjectured by D. Calegari, who also suggested the
problem of trying to calculate stable commutator length exactly in
certain finitely presented groups. I am very grateful to Danny
Calegari for giving generous support and advice during the
preparation of this work.

\section{Stable Commutator length and Quasimorphisms}
\begin{defn}
Let $G$ be a group. A {\em quasimorphism} on $G$ is a function
$$\phi:G \to \R$$
for which there is a constant $D(\phi) \ge 0$ such that for any $a,b
\in G$, we have an inequality
$$|\phi(a) + \phi(b) - \phi(ab)| \le D(\phi)$$
In other words, a quasimorphism is like a homomorphism up to a
bounded error. The least constant $D(\phi)$ with this property is
called the {\em defect} of $\phi$.
\end{defn}

\begin{defn}
A quasimorphism is {\em homogeneous} if it satisfies the
additional property
$$\phi(a^n) = n\phi(a)$$
for all $a \in G$ and $n \in \Z$.
\end{defn}
A homogeneous quasimorphism is a class function by its definition.
Denote the vector space of all homogeneous quasimorphisms on $G$
by $Q(G)$.

\begin{exa}
Let $\widetilde{{\text{Homeo}^+}}(S^1) = \{ f \in
\text{Homeo}^+(\R) | f(x+1) = f(x) + 1 \}$. It's the set
consisting of all possible lifts of elements in
${\text{Homeo}^+}(S^1)$ under the covering projection $\pi : \R
\rightarrow S^1$. We have the central extension:
$$0 \rightarrow \Z \rightarrow
\widetilde{{\text{Homeo}^+}}(S^1) \rightarrow
{\text{Homeo}^+}(S^1) \rightarrow 1$$ where $\Z$ is generated by
the unit translation and $p : \widetilde{{\text{Homeo}^+}}(S^1)
\rightarrow {\text{Homeo}^+}(S^1)$ is the natural projection from
its definition. For $g \in \widetilde{{\text{Homeo}^+}}(S^1)$,
define
$$rot(g) = \lim_{n \to \infty} \frac {g^{n}(0)} n$$
With this definition, $rot$ is a homogeneous quasimorphism with
defect $1$. (See appendix for a proof.)
\end{exa}

We have the following fundamental theorems of C. Bavard, which
state the duality between $scl$ and homogeneous quasimorphisms.
See C. Bavard \cite{ba} or D. Calegari \cite{ca5} for a reference.
\begin{thm}[C. Bavard]\
\begin{enumerate}
\item
Let $G$ be a group. Then for any $a \in [G,G]$, we have an equality
$$scl(a) = \frac{1}{2} \sup_{\phi\in Q(G)/H^1(G;\R)}
\frac{|\phi(a)|}{D(\phi)}$$
\item
There is an exact sequence
$$0 \rightarrow H^1(G;\R) \rightarrow Q(G) \rightarrow H_{b}^2(G;\R)
\rightarrow H^2(G;\R)$$
\item
Let $G$ be a group. Then the canonical map from bounded cohomology
to ordinary cohomology $H_{b}^2(G;\R) \rightarrow H^2(G;\R)$ is
injective if and only if the stable commutator length vanishes on
$[G,G]$.
\end{enumerate}

\end{thm}

Given a group $G$, if $Q(G)$ is small enough, we can use (1) of
Theorem 2.4 to determine $scl$. For the previous example $G =
\widetilde{{\text{Homeo}^+}}(S^1)$, we have $dim_{\R} Q(G) = 1$.
\begin{prop}[J. Barge and \'{E}. Ghys \cite{bar}]\
$\text{rot} : \widetilde{{\text{Homeo}^+}}(S^1) \rightarrow \R$ is
the unique homogeneous quasimorphism which sends the unit
translation to $1$.
\end{prop}
\begin{pf}
Suppose $\tau \in Q(\widetilde{{\text{Homeo}^+}}(S^1))$ is another
such map , then we consider
$$ rot - \tau : \widetilde{{\text{Homeo}^+}}(S^1)
\rightarrow \R$$ which is also a homogeneous quasimorphism, and
since any homogeneous quasimorphism on abelian groups, especially
$\Z\oplus\Z$, must be a homomorphism, we have
$$(rot - \tau)(f_1) = (rot - \tau)(f_2)$$
if $p (f_1)=p (f_2)$. It therefore induces a homogeneous
quasimorphism on $\text{Homeo}^+(S^1)$, denote it still by $rot -
\tau : {\text{Homeo}^+}(S^1) \rightarrow \R$. But
$\text{Homeo}^+(S^1)$ is uniformly perfect, i.e. $cl$ is bounded
(\cite{gh2}), so the induced map is bounded. Since it's
homogeneous, it must be a zero map, i.e.  $rot = \tau$
\end{pf}

By Theorem 2.4 and Proposition 2.5 together, we have $scl(g) =
\frac{1}{2} rot(g)$, for any $g$ in
$\widetilde{{\text{Homeo}^+}}(S^1)$. Suppose $G$ is a subgroup of
${\text{Homeo}}^+(S^1)$ which is uniformly perfect. Let
$\tilde{G}$ be the preimage of G in
$\widetilde{{\text{Homeo}^+}}(S^1)$. Then by the same argument,
$scl(a) = rot(a)/2D$ for any $a$ in $\tilde{G}$, where $D$ is the
defect of $rot$ restricted to $\tilde{G}$.

Our goal therefore is to find the subgroups of
${\text{Homeo}}^+(S^1)$ which are finitely presented and uniformly
perfect, and contain elements with interesting rotation numbers.
As a first example, consider Thompson's well-known group of dyadic
piecewise linear homeomorphisms.

Let $\tilde{G}$ consist of piecewise linear homeomorphisms $f$ of
$\R$ with the following properties:
\begin{enumerate}
\item
For each point $x_i$ of discontinuity of the derivative of $f$
(hereafter a ``break point''), both $x_i$ and $f(x_i)$ are dyadic
rational numbers (i.e. of the form $p2^q, p, q \in \Z$);
\item
The derivatives of the restrictions of $f$ to $(x_i, x_{i+1})$ are
powers of $2$ (i.e. of the form $2^m, m \in \Z$);
\item
$f$ preserves dyadic rational numbers and $f(x+1) = f(x) + 1$.
\end{enumerate}
The elements of $\tilde{G}$ induce piecewise linear homeomorphisms
of $S^1 \simeq \R/\Z$. The collection of these homeomorphisms is
the {\em Thompson group} $G$ (\cite{tho}, \cite{cafp}). Thompson
group is simple, $FP_{\infty}$ and uniformly perfect(\cite{cafp},
\cite{de} and \cite{gh1}), so we have (The defect is still 1. See
appendix.)
$$scl(a) = \frac{1}{2}rot(a),  a \in \tilde{G}$$

About the rotation numbers of elements in the Thompson group, we
have
\begin{thm}[\'{E}. Ghys and V. Sergiescu \cite{gh1}]\
$\text{rot}(a)$ is rational for any $a \in \tilde{G}$.
\end{thm}
So $scl$ takes only rational values on the group $\tilde{G}$. \

\section{Generalized Thompson Groups}
Our definition of generalized Thompson groups is from \cite{st} by
M. Stein. Let $P$ be a multiplicative subgroup of the positive
real numbers and let $A$ be a ${\Z}P$-submodule of the reals with
$PA=A$. Choose a number $l \in A, l>0$. Let $F(l, A, P)$ be the
group of piecewise linear homeomorphisms of $[0, l]$ with finitely
many break points, all in $A$, having slopes only in $P$.
Similarly define $T(l, A, P)$ to be the group of piecewise linear
homeomorphisms of $[0, l]/\{0, l\}$ (the circle formed by
identifying endpoints of the closed interval $[0, l]$) with
finitely many break points in $A$ and slopes in $P$, with the
additional requirement that the homeomorphisms send $A\cap[0, l]$
to itself. In these notations, $T(1, \Z[\frac{1}{2}], \langle
2\rangle)$ is the Thompson group. In our study of generalized
Thompson groups, we always assume that $P$ is generated by the set
of positive integers $\{n_1, n_2, \ldots, n_k \}$ and $A =
\Z[\frac{1}{n_1}, \frac{1}{n_2}, \ldots, \frac{1}{n_k}]$, here
$\{n_1, n_2, \ldots, n_k \}$ forms a basis for $P$. An important
theorem in studying generalized Thompson groups is the following
Bieri-Strebel criterion (See \cite{st} appendix for a proof.).
\begin{thm}[R. Bieri and R. Strebel \cite{bi}]\
Let $a, c, a', c'$ be elements of $A$ with $a<c$ and $a'<c'$. Then
there exists $f$, a piecewise linear homeomorphism of $\R$, with
slopes in $P$ and finitely many break points, all in $A$, mapping
$[a, c]$ onto $[a', c']$ $\Leftrightarrow$ $c'-a'$ is congruent to
$c-a$ modulo $IP*A$.
\end{thm}

Here $IP*A$ is the submodule of $A$ generated by elements of the
form $(1-p)a$, where $a \in A$ and $p \in P$. Let $d = gcd(n_1-1,
\ldots, n_k-1)$ and from now on, we assume that $d = 1$. In this
case, $IP*A = P(d\Z) = P\Z = A$, so the Bieri-Strebel criterion
from Thm 3.1 is vacuously satisfied.

Take an arbitrary $f \in T = T(l, A, P)$ with the assumptions
above. Choose points $a< b, c<d \in [0, l) \cap A$, $f([a, b]) =
[a_1, b_1]$ such that $[c, d] \cap ([a, b]\cup[a_1, b_1]) =
\emptyset$(This can be achieved by taking $[a, b]$ small enough.).
We can find piecewise linear homeomorphisms $g_1$ and $g_2$ with
slopes in $P$ and break points in $A$, sending $[b, c]$ to $[b_1,
c]$ and $[d, a]$ to $[d, a_1]$ respectively. Construct $g \in T$
as follows.
\[ g = \left\{ \begin{array}{ll}
         f & \mbox{if $x \in [a, b]$}\\
         g_1 & \mbox{if $x \in [b, c]$}\\
         id & \mbox{if $x \in [c, d]$}\\
         g_2 & \mbox{if $x \in [d, a]$}\end{array} \right. \]
Now $f = (fg^{-1})g$, both $f\circ g^{-1}$ and $g$ fix some
nonempty open arcs.

Write $f = g_1g_2$ where $g_i$'s fix some open arcs for $i=1,2$. For
each such $g_i$, there exists a rotation through $\theta_i$,
$R_{{\theta}_i} \in T$ such that $supp (R_{\theta_i}\circ g_i\circ
R_{-\theta_{i}}) \subseteq (0, l)$.  ($R_{{\theta}_i} \in T$ if and
only if $\theta_i \in A$ which is dense in $[0, l]$, so we only need
to choose $\theta_i$ in the open arc fixed by $g_i$.). Write $h_i =
R_{\theta_i}\circ g_i\circ R_{-\theta_{i}}$, so $h_i \in F=F(l, A,
P)$ and $h_i$ lies furthermore in the kernel of the following
homomorphism
$$\rho : F \rightarrow P \times P$$
$\varphi \mapsto (\varphi'(0+), \varphi'(l-))$, taking derivatives
at endpoints. Let $B = ker\rho$, so $h_i \in B$.
\begin{thm}[M. Stein \cite{st}]\
$B'$ is simple and $B' = F'$
\end{thm}
\begin{thm}[K. S. Brown]\
$H_*(F) \cong H_*(B)\otimes H_*(P\times P)$.
\end{thm}
(For a proof, see M. Stein \cite{st}.)

Let's further assume that the slope group $P$ has rank 2, i.e. $
P=\langle p, q\rangle, A=\Z[\frac{1}{p}, \frac{1}{q}]$, $p, q \in
\Z_+$ and $d = gcd(p - 1, q - 1) = 1$. $F_{p, q}, T_{p, q}$ are
the corresponding groups. M. Stein, in \cite{st}, computed the
homology groups of $F_{p, q}$ by using its action on a complex.
\begin{thm}[M. Stein]\
$H_1(F_{p,q})$ is a free abelian group with rank $2(d+1)$, where
$d = gcd(p-1, q-1)$.
\end{thm}

By assumption $d = 1$, $rk_{\Z}(H_1(F_{p,q})) = 4$. By Theorem
3.3, $H_1(F_{p,q}) \cong H_1(B)\oplus H_1(P\times P)$ and here
$P\times P \cong \Z\times\Z\times\Z\times\Z$, so $rk_{\Z}(P\times
P) = 4 = rk_{\Z}(H_1(F_{p, q}))$, which implies that $H_1(B) =
B/B'$ is trivial. So $B = B' = F_{p,q}'$.

So $h_i \in B = F_{p,q}'$ can be written as a product of commutators
of $F_{p,q}$. So is $g_i$, which is conjugate to $h_i$ in $T_{p,
q}$. So overall we proved that $T_{p, q}$ is perfect, i.e. $T_{p, q}
= T_{p, q}'$.\\

Let's compute the set $Q(T_{p,q})$. For this purpose, we need the
theorem of D. Calegari \cite{ca3} about the $scl$ of elements in
subgroups of $PL^{+}(I)$.
\begin{thm}[D. Calegari]\
Let $G$ be a subgroup of $\text{PL}^+(I)$. Then the stable
commutator length of every element of $[G, G]$ is zero.
\end{thm}

\begin{lem}
Let $T_{p,q} = T(1, \Z[\frac{1}{p}, \frac{1}{q}], \langle p,
q\rangle)$ and $d = gcd(p-1, q-1) = 1$, then $Q(T_{p, q}) =
\{0\}$.
\end{lem}

\begin{pf}
By the argument above, $f = g_1 g_2$ and $g_i =
R_{-\theta_{i}}\circ h_i\circ R_{\theta_{i}}, i = 1,2$. Since
every homogeneous quasimorphism $\phi$ is a class function, it
follows that $\phi(g_i) = \phi(h_i)$. Since $d = 1$, $h_i \in
F_{p,q}' \subseteq \text{PL}^+(I)$ and Theorem 3.5 and Theorem 2.4
imply $\phi(h_i)=0$, so
 $\phi(g_i)=0$. For any $n$, suppose $f^n = g_{1n}g_{2n}$, then
$\phi(g_{1n}) = \phi(g_{2n}) = 0$
$$ |n\phi(f)| = |\phi(f^n)| = |\phi(g_{1n}g_{2n})| = |\phi(g_{1n}g_{2n}) - \phi(g_{1n})
- \phi(g_{2n})| \leq D(\phi)$$ so
$$|\phi(f)| \leq \frac{D(\phi)}{n}$$
for any $n \in \Z_+$. Let $n \rightarrow +\infty$, we get $\phi(f)
= 0$. By Theorem 2.4, $Q(T_{p, q})/H^1(T_{p, q}; \R) = \{0\}$, but
we have $T_{p, q} = T_{p, q}'$, which implies that $H^1(T_{p, q};
\R)$ is trivial, thus $Q(T_{p, q}) = \{0\}$.
\end{pf}

\begin{thm}
Suppose $gcd(p-1, q-1) = 1$. Let $\tilde{T}_{p,q}$ be the preimage
of $T_{p, q}$ in $\widetilde{{\text{Homeo}^+}}(S^1)$. Then rot is
the unique homogeneous quasimorphism which sends the unit
translation to $1$.
\end{thm}

\begin{pf}
By the same argument in Proposition 2.5 with Lemma 3.6 replacing
the uniformly perfectness.
\end{pf}

By Bavard's Theorem 2.4, if $gcd(p-1, q-1) = 1$, then for all $a$
in $\tilde{T}_{p,q}$ we have an equality $scl(a) = rot(a)/ 2D$. In
the appendix we show that $D = D(rot) = 1$ on the groups
$\tilde{T}_{p,q}$, and therefore this simplifies to $scl(a) =
\frac {1}{2} rot(a)$. So it remains to determine the rotation
numbers of elements in generalized Thompson groups.

Isabelle Liousse (in \cite{li}) has studied the rotation numbers
of elements in generalized Thompson groups, and proved the
following theorem:

\begin{thm}[I. Liousse]\
Let $l \in A=\Z[\frac{1}{p}, \frac{1}{q}]$, $p, q \in \Z_+$
independent and $d = gcd(p-1, q-1)$. Suppose $d$ is prime or $1$,
then $T_{l, (p, q)} = T(l, \Z[\frac{1}{p}, \frac{1}{q}],\langle p,
q\rangle)$ contains elements with irrational rotation numbers. In
particular, if $d = 1$, $T_{l, (p, q)} \cong T_{1, (p, q)}$, for
any $l \in \Z$.
\end{thm}

\begin{exa}[from \cite{li}]\
Let $p = 2, q = 3$, $\tilde{T}_{2, 3}$ satisfies the assumption. We
have the element $f \in T_{2, 3}$:
\[ f = \left\{ \begin{array}{ll}
         \frac{2}{3} x + \frac{2}{3} & \mbox{if $x \in [0, \frac{1}{2}]$}\\
         \frac{4}{3} x - \frac{2}{3} & \mbox{if $x \in [\frac{1}{2}, 1]$}
         \end{array} \right. \]

\begin{figure}[ht]
\centerline{\relabelbox  \epsfxsize 2.0truein
 \epsfbox{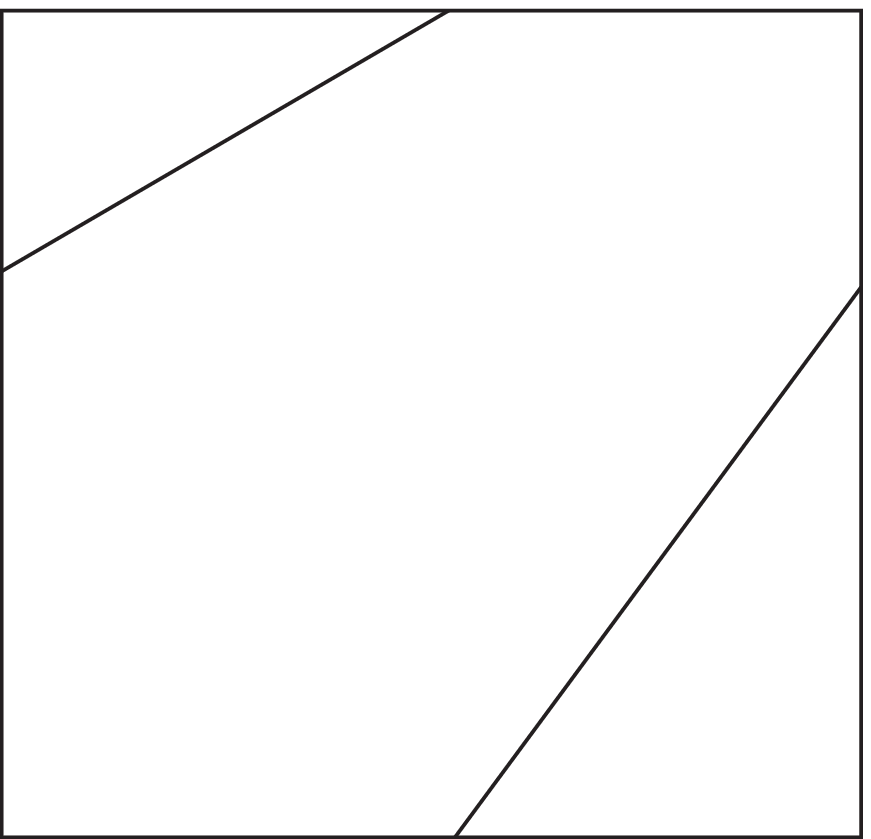}
\extralabel <-155pt,-5pt> {$0$} \extralabel <-155pt,135pt> {$1$}
\extralabel <-75pt,-8pt> {$\frac{1}{2}$} \extralabel <-155pt,95pt>
{$\frac{2}{3}$} \extralabel <0pt,-5pt> {$1$} \extralabel
<-75pt,130pt> {$\frac{1}{2}$} \extralabel <0pt,95pt>
{$\frac{2}{3}$}
\endrelabelbox}
\caption{}\label{ch}
\end{figure}
The transformation $h(x):= 2 - \frac{1}{2^{x-1}}, x \in [0, 1]$
conjugates $f$ to the rotation by $\rho$ where $\rho = \frac{\log
3}{\log 2} -1 \approx 0.58496250072 \cdots$. Thus any lift of $f$
has the rotation number $\frac{\log3}{\log2} + n$, $n \in \Z$. By
the celebrated theorem of Gel'fand-Schneider \cite{la}, the
rotation number $\rho$ is transcendental.

$f$ lies in the class of the homeomorphisms in $PL^+(S^1)$,
satisfying the {\em property D} \cite{li}. I. Liousse showed that
for any $f \in PL^+(S^1)$ with the {\em property D}, there exists
$h \in {\text{Homeo}^+}(S^1)$ such that $h^{-1}fh$ is a rotation,
i.e. in $SO(2)$. Thus we have a measure $\mu$ on $S^1$ that is
invariant under $f$ and furthermore
$$\int \log(Df) d\mu = 0$$
here $Df$ is the function of the derivative of $f$, which is a
piecewise constant function, so to write down the left side of the
equality, we only need to get the measures of the intervals
between break points. This can be obtained by the assumed {\em
property D}, without any knowledge of the conjugate function $h$.
This equality gives a linear equation of the rotation number of
$f$, so in this way, we can get the rotation number.
\end{exa}

\section{Appendix}
In this appendix, we will justify the claims in previous sections
on defects of rotation numbers as homogeneous quasimorphisms. We
will use a lemma by C. Bavard on defect estimation. See C. Bavard
\cite{ba} or D. Calegari \cite{ca5} for a proof.

\begin{lem}[C. Bavard]\
Let $\phi$ be a homogeneous quasimorphism on $G$. Then there is an
equality
$$\sup_{a, b \in G} |\phi([a, b])| = D(\phi)$$
\end{lem}

\begin{prop}
$rot: \widetilde{{\text{Homeo}^+}}(S^1) \to \R$ is a homogeneous
quasimorphism with defect $1$.
\end{prop}

\begin{pf}
Refer to \cite{ka} for basic properties of rotation numbers.\\

(1)\ Let $f, g \in \widetilde{{\text{Homeo}^+}}(S^1)$. Without
loss of generality, we can assume that $0\leq f(0), g(0) < 1$. So
$0\leq f\circ g(0)< 2$. And $0 \leq rot(f) \leq 1$, $0 \leq rot(g)
\leq 1$, and $0 \leq rot(f\circ g) \leq 2$. Thus we have $|
rot(f\circ g) - rot(f) - rot(g)| \leq 2$, $rot$ is a
quasimorphism. That $rot$ is homogeneous
is clear from its definition.\\

(2)\ We show that $D(rot) = 1$ by using Lemma 4.1. \\
Take any $f, g \in \widetilde{{\text{Homeo}^+}}(S^1)$ and we are
going to compute $rot([f, g])$. We can still assume that $0\leq
f(0), g(0) < 1$. Suppose $0 \leq g(0) \leq f(0) < 1$, then we have
by using that $f, g$ are increasing functions:
$$g(f(0)) < g(1) = g(0) + 1 \leq f(0) + 1 \leq f(g(0)) +
1$$ so
$$f(g(0)) - g(f(0)) > -1$$
Then we have two cases:\\
(i) If we also have $f(g(0)) - g(f(0)) \leq 1$, then
$$-1 \leq f(g(0)) - g(f(0)) \leq 1$$
$$g(f(-1)) = g(f(0)) - 1 \leq f(g(0)) \leq g(f(0)) + 1 = g(f(1))$$
which implies
$$-1 \leq f^{-1}g^{-1}fg(0) \leq 1$$
so
$$|rot([f, g])| \leq 1$$\\
(ii) If instead we have $f(g(0)) > g(f(0)) + 1$, then $g(f(0)) <
f(g(0)) - 1 = f(g(0) - 1) < f(0)$. Consider $H(x) =
f^{-1}g^{-1}fg(x) - 1 - x$, for $x \in [0, 1]$. $H(0) =
f^{-1}g^{-1}fg(0) - 1 > 0$ by assumption.
{\setlength\arraycolsep{2pt}
 \begin{eqnarray}
 H(f(0)) & = & f^{-1}g^{-1}fg(f(0)) -
1 - f(0) \nonumber  \\
 & < & f^{-1}g^{-1}f(f(0)) -
1 - f(0) \nonumber  \\
 & = & f^{-1}g^{-1}f^2(0) -
1 - f(0) \nonumber  \\
\nonumber
\end{eqnarray}}
We want to show that $H(f(0)) < 0$, which can be deduced from the
inequality below
$$f^{-1}g^{-1}f^2(0) < 1 + f(0)$$
which is equivalent to
$$f^2(0) < g(f^2(0)) + 1$$
This is always true since $x < g(x) + 1$, for any $x \in \R$.\\
So we have $H(0) > 0$ and $H(f(0)) < 0$, here $0 < f(0) < 1$. There
must be a point $y \in (0, f(0))$ such that\\
$$H(y) = f^{-1}g^{-1}fg(y) - 1 - y = 0$$
that is
$$f^{-1}g^{-1}fg(y) = 1 + y$$
So
$$rot[f,g] = \lim_{n \to \infty} \frac{[f, g]^n(y) - y}{n} = 1$$\\
\\
The proof for the case $0 \leq f(0) \leq g(0) < 1$ is the same. Put
all together and we get $D(rot) = 1$ by C. Bavard's Lemma 4.1.
\end{pf}

\begin{prop}
Let $T = T(1, A, P)$, where $A = \Z[\frac{1}{n_1}, \frac{1}{n_2},
\ldots, \frac{1}{n_k}]$, $P = \langle n_1, n_2, \ldots, n_k
\rangle$, here $n_i$'s are independent. Suppose $d = gcd(n_1 - 1,
\ldots, n_k - 1) = 1$, then $T$ is dense in $\text{Homeo}^+(S^1)$
with {\em $C^0$} topology.
\end{prop}

\begin{pf}
Take an arbitrary $f$ in $\text{Homeo}^+(S^1)$ and any $\epsilon >
0$. $f$ is uniformly continuous, so we can choose $0= x_0 < x_1 <
\ldots < x_l <1$, $x_i$'s are in $A$ that is dense in $S^1 = [0,
1]/\{0=1\}$, such that $|f(x_i) - f(x_{i-1})| <
\frac{\epsilon}{3}$. Since $A$ is dense in $[0, 1]$, we can find
$y_0 < y_1 < \ldots < y_l <1$, $y_i$'s are in $A$ and $|y_i -
f(x_i)| < \frac{\epsilon}{3}$. By the Bieri-Strebel
criterion(Theorem 3.1), there exists $g \in T$ such that $g(x_i) =
y_i$. And from the choice of $x_i$'s and $y_i$'s, it's easy to see
that $\|f - g\|_{C^0} < \epsilon$.
\end{pf}

Thus we also have that $\tilde{T} \subseteq
\widetilde{{\text{Homeo}^+}}(S^1)$ is dense. On the other hand,
the map $rot$, thought of as a function from
$\widetilde{{\text{Homeo}^+}}(S^1)$ to $\R$, is continuous in the
{\em $C^0$} topology \cite{ka}. So we have that the defect of
rotation number $D(rot)$, restricted to $\tilde{T}$, is also 1.

\end{document}